\documentclass[reqno]{amsart}
\usepackage{amsthm}
\usepackage{amsmath}
\usepackage{amssymb}
\usepackage{amsrefs}
\usepackage{enumitem}
\usepackage{crossreftools}
\usepackage{graphicx}
\usepackage{rotating}

\newcommand{\ve}{\varepsilon}

\newcommand{\om}{\omega}
\newcommand{\be}{\begin{equation}}
\newcommand{\ee}{\end{equation}}
\newcommand{\ba}{\begin{align}}
\newcommand{\ea}{\end{align}}

\newcommand\TS{\rule{0pt}{2.6ex}}
 
\DeclareMathOperator{\Res}{Res}

\DeclareMathOperator{\GA}{GA}

\newtheorem{theorem}{Theorem}
\newtheorem*{theorem*}{Theorem}

\title[Integer group determinants]{The integer group determinants for the abelian groups of order 27}

\author[M. Mossinghoff]{Michael J.\ Mossinghoff}
\address{Center for Communications Research\\
Princeton, NJ USA}
\email{m.mossinghoff@ccr-princeton.org}

\author[C. Pinner]{Chris Pinner}
\address{Department of Mathematics\\
Kansas State University\\
Manhattan, KS 66506, USA}
\email{pinner@math.ksu.edu}


\date{\today}

\begin{document}

\begin{abstract}
We completely determine the integer group determinants for $\mathbb Z_3 \times \mathbb Z_9$ by resolving the allowable multiples of~$3$.
This completes the Taussky-Todd integer group determinant problem for the abelian groups of order $27$.
\end{abstract}

\maketitle

\section{Introduction}\label{secIntro}
At the April 1977 meeting of the AMS  in Hayward, California, Olga Taussky-Todd \cite{TausskyTodd} asked whether one could determine the values taken by the group determinant of a group $G$ when the entries are all integers.
We shall think of the group determinant as a function on the group ring $\mathbb{Z}[G]$,
\[
D_G\Biggl(\sum_{g\in G} a_g g\Biggr)= \det \left( a_{g_ig_{j}^{-1}}\right),
\]
and the Taussky-Todd integer group determinant problem for $G$ is to determine
\[
\mathcal{S}(G)=\Biggl\{ D_G\Bigl( \sum_{g\in G} a_g g\Bigr)\; :\; a_{g}\in \mathbb Z \Biggr\}.
\]
In the abelian case, the group representations are just the linear characters, and the classic Frobenius \cite{Frob} factorisation of the group determinant (see for example Conrad \cite{Conrad} or Johnson \cite{Johnson}) takes a particularly simple form, as was already known to Dedekind.
We let $\mathbb{Z}_n$ denote the cyclic group of order $n$.
If
\[
G=\mathbb Z_{n_1}\times \cdots \times \mathbb Z_{n_k},
\]
so that
\[
\mathbb Z[G]= \mathbb Z [x_1,\ldots , x_k]/\langle x_1^{n_1}-1,\ldots ,x_k^{n_k}-1\rangle,
\]
then
\[
D_G(f(x_1,\ldots ,x_k)= \prod_{j_1=0}^{n_1-1}\cdots \prod_{j_k=0}^{n_k-1} f(\om_{n_1}^{j_1},\ldots ,\om_{n_k}^{j_k}),\quad \om_{\ell}:=e^{2\pi i/\ell}.
\]
Rapid progress was made \cite{Laquer,Newman1} for the cyclic groups $\mathbb Z_p$ and $\mathbb Z_{2p}$,  where straightforward necessary divisibility conditions turn out to be sufficient as well, and likewise \cite{Mahoney2,dihedral} for the dihedral groups of order $2p$ and $4p$.
However, Newman \cite{Newman2} showed that this was not true even for $\mathbb Z_{p^2}$, except for $p= 3$.
In \cite{MP23} we explored some of the complexity encountered with $\mathbb Z_{p^k}$ when $k\geq 2$: there we obtained complete descriptions for $\mathbb Z_{25}$ and $\mathbb Z_{27}$, but even these were not straightforward.

For $\mathbb Z_{27}$, we needed to divide the primes $p\equiv 1\bmod 3$ into two classes.
These primes split in $\mathbb Q(\om_3)$, and we may write
\be \label{Type1def} p=N_1( A+3B(x-1)), \quad A,B\in \mathbb Z, \ee
where $N_1(f(x))=f(\om_3)f(\om_3^2)=\Res(x^2+x+1,f(x))$ is the $3$-norm from $\mathbb Q(\om_3)$ to $\mathbb Q$.
Here and throughout $\Res(f(x), g(x))$ denotes the resultant of the polynomials $f(x)$ and $g(x)$.
We call the prime $p$ \textit{Type~1} if $3\nmid B$ and \textit{Type~2} if $3\mid B$.
The Type~2 primes are the OEIS \cite{OEIS} sequence A014753.
The Type~2  primes $p\equiv 1 \bmod 9$ split completely in $\mathbb Q(\om_9)$, and we may write
\be \label{Type2def} p=N_2(1+ (x-1)^7h(x)),\quad h(x)\in \mathbb Z[x],  \ee
where $N_2(f(x))=\Res(x^6+x^3+1,f(x))$ is the $9$-norm from $\mathbb Q(\om_9)$ to $\mathbb Q$.
We call these primes \textit{Type~3} if $3\nmid h(1)$ and \textit{Type~4} if $3\mid h(1)$.
We describe an algorithm that classifies an element in $\mathbb Z[\om_3]$ or $\mathbb Z[\om_9]$ to the form \eqref{Type1def} or \eqref{Type2def} in Section~\ref{secAlg}.

In 2023 we established the following statement.

\begin{theorem*}[\cite{MP23}, Theorem 3.1]
The integer group determinants for $\mathbb Z_{27}$ are the integers coprime to $3$, all multiples of $3^6$, the integers $3^4k$ with $k$ divisible by a Type~1 prime $p\equiv 1\bmod 9$, and the integers $3^5k$ with $k$ divisible by either a Type~1 prime $p\equiv 1\bmod 3$ or a Type~3 prime $p\equiv 1\bmod 9$.
\end{theorem*}

A complete description of $\mathcal{S}(G)$ has been obtained for all groups of order less than $21$: see \cite{smallgps} for $|G|\leq 14$, \cite{bishnu1} for $\mathcal{S}(\mathbb Z_{15})$, \cite{YYCirc16,YYC8C2,YYC24,YYC42,YYC422} for the abelian groups of order $16$, \cite{dihedral,PaudelPinnerZnxH,bishnu4,PPQ16,YYNA16,YY16} for the non-abelian groups of this order, and \cite{PaudelPinner25} for the groups of order $18$.
The non-abelian group of order $21$, along with $\GA(1,5)$ of order $20$, is dealt with in \cite{bishnu2}. The remaining groups of order 20 are tackled in \cite{Radolovic}.
The abelian groups of order $24$ are addressed in \cite{Panraksa24}, and $S_4$ is treated in \cite{S4}.
This leaves only twelve groups of order less than $27$ unresolved: the cyclic group $\mathbb{Z}_{21}$ and eleven non-abelian groups of order $24$. 

For groups of the form $\mathbb Z_p^n$, the values coprime to $p$ are straightforward: as shown in \cite{dilum1}, they are the solutions to $x^{p-1}\equiv 1 \bmod {p^n}$.
The multiples of $p$, on the other hand, are typically much more complicated to describe, other than $\mathbb Z_5^2$ \cite{Panraksa25} (see also \cite{Panraksa49} for $\mathbb Z_7^2$).
In \cite{MP26} we derived bounds on the smallest  power of $p$ dividing a group determinant for $G=\mathbb Z_{p}^n$ and other products of abelian $p$-groups, and on the smallest integer $k$ with $p^k\mathbb Z \subset \mathcal{S}(G)$ for such $p$-groups.
Recently, Panraksa \cite{Pan2} showed that there is subtlety in the exponent $k,$ even for groups of the form $G=\mathbb Z_p^2$.
In \cite{MP26}, we obtained a complete description of the integer group determinants for $G=\mathbb Z_3^3$.

\begin{theorem*}[\cite{MP26}, Theorem 5.1]
For $G=\mathbb Z_3^3$, the integer group determinants coprime to $3$ are the integers $D\equiv \pm 1 \bmod 27$.
The set of group determinants divisible by $3$ consists of
\begin{enumerate}[label=(\alph*)]
\item all multiples of $3^{20}$,
\item the integers $3^{19}m$ with $m$ divisible by a prime $p\equiv 4$ or $7 \bmod 9$, or the square of a prime $p^2$ with $p\equiv 2$ or $5 \bmod 9$,
\item the integers $3^{18}m$ where either $m\equiv \pm 1$ or $\pm 2 \bmod 9$, or $m\equiv \pm 4 \bmod 9$ with $m$ divisible by at least one of the following: a prime $p\equiv 4$ or $7 \bmod 9$, a Type~1 prime $p\equiv 1 \bmod 9$, or $p^2$ with $p\equiv 2$ or $5 \bmod 9$.
\end{enumerate}
\end{theorem*}

When the products contain $\mathbb Z_{p^t}$ with $t\geq 2$, we no longer have a simple description of the values coprime to $p$.
In \cite{MP26}, we illustrated this by determining the group determinants coprime to $3$ for the group $\mathbb Z_3 \times \mathbb Z_9$.
There we also obtained sharp $3$-adic bounds for $\mathcal{S}(\mathbb{Z}_3\times\mathbb{Z}_9)$.
We summarize these results.

\begin{theorem*}[\cite{MP26}, Theorems 3.2(iii) and 3.4]
For the group $\mathbb Z_3 \times \mathbb Z_9$, the integer group determinants coprime to $3$ are the integers $D\equiv \pm 1\bmod 27$, together with the integers $D\equiv \pm 8$ or $\pm 10\bmod 27$ of the form $D=A(B^2-BC+C^2)$, for some integers $A$, $B$, $C$ satisfying
\[
A\equiv \pm 2 \text{ or } \pm 4 \bmod 9,\quad B\equiv A \bmod 9,\quad C\equiv A-B\bmod 27.
\]
For the group determinants $n\in\mathcal{S}(\mathbb{Z}_3\times\mathbb{Z}_9)$ divisible by $3$, we have $3^9\mid n$, and furthermore $3^{11}\mathbb{Z} \subset \mathcal{S}(\mathbb{Z}_3\times\mathbb{Z}_9)$.
\end{theorem*}

In order to complete the description of $\mathcal{S}(\mathbb{Z}_3\times\mathbb{Z}_9)$, it only remains to characterize the multiples of $3^9$ and $3^{10}$ in this set.
The purpose of this article is to determine these by proving the following theorem.

\begin{theorem}\label{main}
For the group $\mathbb{Z}_3\times\mathbb{Z}_9$, the group determinants divisible by $3$ are
\begin{enumerate}[label=(\alph*)]
\item\label{main11} all multiples of $3^{11}$,
\item\label{main10} integers $3^{10}m$ where $m$ is divisible by a Type~1 prime $p\equiv 1 \bmod 9$, and
\item\label{main9} integers $3^9m$ where $m$ is of the form
\begin{enumerate}[label=(\roman*)]
\item\label{main9i} $m\equiv \pm 1 \bmod 9$, or
\item\label{main9ii} $m\equiv \pm 4 \bmod 9$, or
\item\label{main9iii} $m\equiv \pm 7 \bmod 9$ where $m$ is divisible by a Type~1 prime $p\equiv 1 \bmod 3$, or a Type~3 prime $p\equiv 1 \bmod 9$.
\end{enumerate}
\end{enumerate}
\end{theorem}

It is interesting that the same classification of primes as Type~1, 2, 3 or~4 that occurred in the characterizations for $\mathcal{S}(\mathbb{Z}_{27})$ and $\mathcal{S}(\mathbb{Z}_3^3)$ is exactly what is required for that of $\mathbb Z_3 \times \mathbb Z_9$.
The  Type~1 primes less than $2500$ are shown below.
\begin{quote}
7, 13, 19, 31, 37, 43, 79, 97, 109, 127, 139, 157, 163, 181, 199,
211, 223, 229, 241, 277, 283, 313, 331, 337, 349, 373, 379, 397, 409,
421, 433, 457, 463, 487, 541, 571, 601, 607, 631, 673, 691, 709, 733,
739, 751, 769, 811, 823, 829, 859, 877, 883, 907, 937, 1033, 1039,
1051, 1063, 1069, 1087, 1123, 1129, 1153, 1171, 1201, 1213, 1231,
1237, 1279, 1291, 1297, 1327, 1381, 1423, 1429, 1447, 1453, 1459,
1471, 1483, 1489, 1567, 1579, 1627, 1657, 1663, 1693, 1699, 1723,
1741, 1747, 1753, 1777, 1789, 1801, 1831, 1873, 1951, 1987, 1993,
1999, 2011, 2017, 2053, 2083, 2089, 2113, 2137, 2143, 2161, 2239,
2281, 2293, 2311, 2347, 2371, 2377, 2437, 2467, 2473.
\end{quote}
The Type~3 primes less than $5000$ that are $1\bmod 9$ are
\begin{quote}
73, 271, 307, 523, 577, 613, 757, 919, 1009, 1531, 1621, 1783, 2179, 2287, 2971, 3079, 3187, 3529, 3853, 3889, 4177, 4339, 4951, 4969,
\end{quote}
and the Type~4 primes in the same range and congruence class are
\begin{quote}
991, 1117, 1549, 2251, 2269, 2341, 3907, 4483, 4933.
 \end{quote}
See \cite{GithubData} for more extensive lists.
We established in \cite[Lemma 4.9]{MP23} that asymptotically $2/3$ of the primes $p\equiv 1$ mod $3$ are Type~1.
The Type~2 primes can also be characterized as the primes $p\equiv 1\bmod 3$ where $3$ is a cube mod $p$.
The Type~4 primes appear in \cite{Ozaki}, where they are defined as Type~2 primes $p\equiv 1 \bmod 9$ with the additional property that $(x+1)(x^{-1}+1)$ is also a cube mod $p$ for any $x$ of order $9 \bmod p$.

Theorem~\ref{main} completes the determination of $\mathcal{S}(G)$ for the abelian groups of order $27$.
There are two non-abelian groups of this size.
For one of these, the Heisenberg group $H_3$ of order $27$, in 2024 \cite[Theorem 3.6]{Heisenberg} we proved that
\[
\mathcal{S}(H_3)=\{m\equiv \pm 1 \bmod 27\} \cup 3^{12}\mathbb Z.
\]
We thus have the proper inclusions
\[
\mathcal{S}(\mathbb Z_{3}^3 ) \subset \mathcal{S}(H_3) \subset \mathcal{S}(\mathbb Z_3\times \mathbb Z_9) \subset \mathcal{S}(\mathbb Z_{27}).
\]
This leaves only one group of order $27$ to be investigated: $\mathbb Z_9 \rtimes \mathbb Z_3$.

This article is organized in the following way. Section~\ref{secAlg} develops algorithms for classifying primes as Type~1, 2, 3 or~4, and discusses what it means for an integer to be a 3-norm or 9-norm of one of these types.
Section~\ref{secAchieving} shows that all of the values listed in parts~\ref{main10} and~\ref{main9} of Theorem~\ref{main} are achieved as integer group determinants of $\mathbb{Z}_3\times\mathbb{Z}_9$.
Section~\ref{secNoMore10}
shows that the values in part~\ref{main10} are the only allowable values $D$ of the integer group determinants for $\mathbb{Z}_3\times\mathbb{Z}_9$ with $3^{10}\parallel D$.
Finally,  Section~\ref{secNoMore9} establishes the corresponding statement for part~\ref{main9} of Theorem~\ref{main}, using a computation that employs the algorithms developed in Section~\ref{secAlg}.

\section{Algorithms to determine Type~1, 2, 3 or 4}\label{secAlg}

\subsection{3-norms} 
In \cite[Lemma 4.1]{MP23}, we showed that an element $\alpha$ of $\mathbb Z[\om_3]$ with $3\nmid N_1(\alpha)$ can be written in the form
\be \label{basicform}  u\alpha = A + 3B(x-1), \quad A,B\in \mathbb Z, \ee
with $u$ a unit and $x$ a primitive cube root of unity.
There we classified the algebraic integer $\alpha$ as Type~1 if $3\nmid B$ and Type~2 if $3\mid B$, and we say the $3$-norm $m=N_1(\alpha)$ is a \textit{Type~1} or \textit{Type~2 representation} for $\alpha$ depending on whether $\alpha$ is Type~1 or~2.
From \cite[Lemma 4.2]{MP23}, an Eisenstein integer $\alpha$ cannot be both Type~1 and Type~2, and conjugates must be of the same type.
This aligns with our prior definition of Type~1 and Type~2 primes from \eqref{Type1def}, as the type of $p\equiv 1\bmod 3$ matches that of $\alpha$ when $p=N_1(\alpha)$.
In \cite[Lemmas 4.7 and 4.8]{MP23}, it was shown that if $m=N_1(\alpha)$, $3\nmid m$, and $m$ contains no Type~1 primes, then $\alpha$ must be Type~2.
(It is Type~1 if $m$ contains exactly one Type~1 prime; both types are possible otherwise depending on the conjugates used to build $\alpha$.)

As shown in \cite[Lemma 4.1]{MP23}, the condition \eqref{basicform} is easy to check.
For a primitive cube root of unity $x$, we can first use reduction mod $x^2+x+1$ to write 
\be \label{Type1form} \alpha = a + b(x-1) \ee 
with $3\nmid N_1(\alpha)$, guaranteeing that $3\nmid a$.
If $3\nmid b$ then we can take 
\[
x\alpha =(a-3b)+ (a-2b)(x-1)  \quad \text{ or } \quad  x^2\alpha =(3b-2a)+(b-a)(x-1).
\]
In particular if $\alpha$  is Type~2, for example if $m=N_1(\alpha)$ contains no Type~1 primes $p\equiv 1$ mod 3,  then we obtain
\begin{equation} \label{Step2} 
\begin{split}
3\mid b  & \implies 9\mid b,\\
b\equiv a \bmod 3 & \implies a\equiv b \bmod 9,\\
b\equiv -a \bmod 3 & \implies a\equiv 2b \bmod 9.
\end{split}
\end{equation}

\subsection{9-norms}\label{subsec9norms}
For a primitive ninth root of unity $x,$ an algorithm was given in \cite[Lemma 4.1]{MP23} to write an $\alpha$ in $\mathbb Z [x]$ with $3\nmid N_2(\alpha)$ in the form
\be \label{Type1exp} u\alpha = 1 + (x-1)^5h(x), \quad h(x)\in \mathbb Z[x],\quad \text{ $u$ a unit}.\ee
We called this a Type~1 expansion if $3\nmid h(1)$ or Type~2 if $3\mid h(1)$.
Further, a Type~2 expansion can be written
\be \label{Type3exp} u'\alpha = 1 + (x-1)^7t(x), \quad t(x)\in \mathbb Z[x],\quad \text{ $u'$ a unit}, \ee
and we called this Type~3 if $3\nmid t(1)$ and Type~4 if $3\mid t(1)$.
We shall call $m=N_2(\alpha)$ a Type~1, 2, 3 or~4 representation of $m$ according to the classification of $\alpha$ into these types.
Here, an $\alpha$ cannot be both Type~1 and~2, or both~3 and~4, and conjugates always have the same type.
A prime $p\equiv 1 \bmod 9$  is Type~1, 3 or~4 as $p=N_2(\alpha)$ with $\alpha$ of that type.
(This is consistent with the $3$-norm definition of Type~1 or 2 for these primes \cite[Lemma 4.4]{MP23}.)
If $m$ is a $9$-norm $m=N_2(\alpha)$ with $3\nmid m$,  and $m$ contains no Type~1 primes, then $\alpha$ must be Type~2, if further $m$ contains no Type~3 primes, then $\alpha$ must be Type~4 (see \cite[Lemmas 4.7 and 4.8]{MP23}).

We describe the algorithm from \cite[Lemma 4.1]{MP23}  to write an $\alpha$ in $\mathbb Z[x]$ with $3\nmid N_2(\alpha)$  in the form \eqref{Type1exp}, and in the form \eqref{Type3exp} for Type~2.
We set $\pi=x-1$.
Using
\[
3= -\pi^6 +\pi^9\bmod  {\pi^{10}},\quad   3\pi^k =0\bmod {\pi^8} \text{ for } k\geq 2,
\]
we first  write 
\be \label{standardform} \pm \alpha = 1 + \sum_{j=1}^7 \ve_j \pi^j \bmod {\pi^8},\quad \ve_j\in \{ 0,-1,1\}. \ee
We note the following eight units and their expansions:
\begin{align*}
u_1 &= x=1+\pi,\\
v_1 &= -(1+x)=1-\pi +\pi^6 \bmod {\pi^8},\\
u_2 &= x^2(1+x)^2=1+\pi^2+\pi^4-\pi^6-\pi^7 \bmod {\pi^8}, \\
v_2 &= -x(1+x)=1-\pi^2+\pi^6+\pi^7 \bmod {\pi^8},   \\
u_3 &=   x^3=1+\pi^3-\pi^7 \bmod {\pi^8},\\
v_3 &= x^6=1-\pi^3+\pi^6+\pi^7 \bmod {\pi^8}, \\
u_4 &= -x^6(1+x)^2(1+x^4)=1+\pi^4+\pi^5-\pi^6 -\pi^7 \bmod{ \pi^8},\\
v_4 &= x^{12}(1+x)^4(1+x^4)^2=1-\pi^4-\pi^5+\pi^6+\pi^7 \bmod {\pi^8}.
\end{align*}
Multiplying \eqref{standardform}  by $v_j$ or $u_j$ (and reducing the resulting coefficients mod $3$) we can successively remove any  $\pm \pi^j$ term, $j=1,\ldots,4,$ until we reach
\[
u\alpha = 1\pm \pi^5+ \ve'_6\pi^6+\ve_7'\pi^7 \bmod \pi^8 \;\; \quad \mathrm{(Type~1)},
\]
or
\[
u\alpha = 1+\ve'_{6}\pi^6+\ve'_7\pi^7  \bmod {\pi^8} \;\; \quad \mathrm{(Type~2)}.
\]
For those of Type~2, we can continue with the units
\begin{align*}
u_5 &= x^6(1+x)^6=1+\pi^6 \bmod {\pi^8},\\
v_5 &= -x^3(1+x)^3=1-\pi^6 \bmod {\pi^8},
\end{align*}
to  remove any $\pm \pi^6$ term and obtain
\[
u\alpha = 1\pm  \pi^7 \bmod {\pi^8}  \;\;\quad  \mathrm{(Type~3)},
\]
or
\[
u\alpha = 1 \bmod {\pi^8} \;\; \quad \mathrm{(Type~4).}
\]
Using this algorithm on the $3^5$ choices of $(\ve_1,\ldots ,\ve_5)$  with $\ve_i\in \{0,-1,1\}$ produced a database of $162$ sequences giving a Type~1 expansion, and $81$ giving a Type~2 expansion.
The $9\cdot 81=729$ choices of $(\ve_1,\ldots ,\ve_7)$ producing a Type~2 were then tested to produce a database of $486$ sequences of length $7$ that gave Type~3, and $243$ that produced Type~4.
All four complete lists are available at \cite{GithubData}.

\section{Achieving the values}\label{secAchieving}
For convenience we treat the group as $\mathbb Z_9\times \mathbb Z_3$ and write
\[
\mathbb Z[\mathbb Z_9\times \mathbb Z_3]=\mathbb Z[x,y]/\langle x^9-1,y^3-1 \rangle.
\]
As noted in the introduction, from the work in \cite{MP26} we only need to consider the integer group determinants $D$ with $3^{10}\parallel D$ and $3^{9}\parallel D$.
We show here that all of the values described in parts~\ref{main10} and~\ref{main9} of Theorem~\ref{main} are attained as values of $D$.

\subsection{Achieving the $3^{10}\parallel D$}\label{subsecAchieve10}
Suppose that $m$ is a $9$-norm divisible by a Type~1 prime $p\equiv 1\bmod 9$, then by \cite[Lemma 4.7]{MP23} it has a Type~1 representation, and by \cite[Lemma 4.3]{MP23} we may write
\[
3m = N_2(1-x+3\delta(x)+3(x-1)^2t(x)), \quad \delta(x)=x \text{ or } -1, \quad t(x)\in \mathbb Z [x].
\]
Since
\[
N_2(1+2x+3(x-1)^2t(x)) = N_2(-x( 1+2/x+3(x^8-1)^2t(x^8))),
\]
we can assume that
\[
3m = N_2(-2-x+3(x-1)^2t(x)), \quad t(x)\in \mathbb Z [x].
\]
If $3\mid t(1)$ then we can replace $(x-1)^2t(x)$ by $(x-1)^3t_1(x)$, and since $x$ is a primitive ninth root of unity we have $(x-1)^3=(x^3-1)-3x(x-1)=(x^3-1)(1-x(x-1)(x^6-1))$, so we can write
\be \label{type1} 3m=N_2(\alpha_1+3(x^3-1)t_2(x)), \quad \alpha_1=-2-x \ee
for some $t_2(x)$ in $\mathbb Z[x]$.
Similarly, if $t(1)\equiv 1\bmod 3$ then we can write
\be \label{type3} 3m=N_2(\alpha_2+3(x^3-1)t_2(x)),\quad \alpha_2=-2-x+3(1-x)^2,  \ee
and if $t\equiv -1\bmod 3$ then
\be \label{type5} 3m=N_2(\alpha_3+3(x^3-1)t_2(x)),\quad \alpha_3=2+x+3(1-x)^2.  \ee
We write 
\[
F(x,y)=f_1(x) + f_2(x)y+f_3(x)y^2
\]
and take
\[
G(x,y) = F(x,y) + (y^2+y+1)(x^3-1)t_2(x) +\ell \left(\frac{y^3-1}{y-1}\right)\left(\frac{x^9-1}{x-1}\right),
\]
so that  $G(1,1)=F(1,1)+27\ell$ and  $G(x,1)=F(x,1)+3(x^3-1)t_2(x)$, for $x$ a primitive ninth root of unity, but $G(x,y)=F(x,y)$ for the other pairs of ninth roots of unity $x$ and cube  roots of unity $y$.
We take the following choices of $f_1(x)$, $f_2(x)$, $f_3(x)$ for \eqref{type1}, \eqref{type3} and \eqref{type5} respectively.
\begin{center}
\begin{tabular}{ccc}
$f_1(x)$ & $f_2(x)$ & $f_3(x)$ \\ \hline\TS
$x^7+x^6$ & $x^8+x^3+x^2$  & $x^6+x^5+x^4+x^3$\\
$x^7+x^4+x^3+x^2-2x$  & $x^4+x^2-x+1$  & $2x^7+x^6+x^4+x^2-x+1$\\
$-x^6+x^4+x^2-x+1$ & 2$x^7+x^4+x^2-x+1$ & $x^7+x^4-x^3+x^2 +2$
\end{tabular}
\end{center}
These choices produce the following respective values for $F(x,1)$:
\[
\alpha_1+(x^2+x+2)\Phi_9(x), \quad
\alpha_2+(3x+1)\Phi_9(x), \quad
\alpha_3 +(3x-1)\Phi_9(x).
\]
In all three cases, we have $G(1,1)=9(1+3\ell)$, $\Res(x^2+x+1,F(x,1))=3$, and $\Res(x^6+x^3+1,G(x,1))=3m$,
while for a primitive cube root of unity $\omega$ one may check that $\Res(x^9-1,F(x,\om)F(x,\om^2))=3^6$.
Hence using $\pm F$ we obtain all $\pm (3\ell+1)3^{10}m$ with $m$ a  $9$-norm of Type~1.
From \cite[Lemmas 4.7 and 4.8]{MP23}, a Type~1 $9$-norm $m$ must contain a Type~1 prime $p\equiv 1\bmod 9$, so we achieve all these values if we restrict ourselves to taking $m$ to be a Type~1 prime $p\equiv 1\bmod 9$.

\subsection{Achieving the $3^{9}\parallel D$}\label{subsecAchieve9}
For~\ref{main9i}, we obtain $\pm 3^9(1+9\ell)$ from $\pm F(x,y)$ of the form
\[
F(x,y) =x^3+x+1  + \ell \left(\frac{y^3-1}{y-1}\right)\left(\frac{x^9-1}{x-1}\right).
\]
For~\ref{main9ii}, we obtain $\pm 3^9(4+9\ell)$ from $\pm F(x,y)$ with
\[
F(x,y)=1+x+x^2+x^3+x^7+x^8 +y(x+x^2+x^8)+y^2(1+x+x^2)+  \ell \left(\frac{y^3-1}{y-1}\right)\left(\frac{x^9-1}{x-1}\right).
\]
For \ref{main9iii}, from \cite[Lemma 4.3]{MP23} any $3$-norm $m$ can be written in the form
\begin{equation}\label{standardreps}
3m = N_1(\delta(1-x)+9a+9b(x-1)), 
\end{equation}
with $\delta=1$, $2$ or $4$ as $m=1$, $4$ or $7 \bmod 9$.
This corresponds to a Type~1 representation for $m$ if $3\nmid a$.
Notice that
\begin{align*} 
N_1(\delta(1-x)+9a+9b(x-1)) &= N_1(-x(\delta(1-1/x)+9a+9b(1/x-1)))\\
 &= N_1(\delta(1-x)-9a+9(b-a)(x-1)),
\end{align*}
so that an  $m$ with a Type~1 representation will have representations \eqref{standardreps} with both $a\equiv \pm 1\bmod 3$.
Note a $3$-norm must be a product of primes $1$, $4$ or $7\bmod 9$ and squares of primes $2$, $5$ or $8\bmod 9$.
From \cite[Lemmas 4.7 and 4.8]{MP23}, a $3$-norm having a Type~1 representation must contain a Type~1 prime $1$, $4$ or $7\bmod 9$.
We show that we can achieve $\pm 3^9mk_1$ for any $3$-norm $m$ with a Type~1 representation and any $k_1$ with $mk_1\equiv 7\bmod 9$. 
Restricting $m$ to a Type~1 prime, we thus get all $k\equiv \pm 7\bmod 9$ with $k$ divisible by a Type~1 prime $p\equiv 1$, $4$ or $7\bmod 9$.

To obtain the $m\equiv 7\bmod 9$ with a Type~1 representation, we set
\[
F(x,y)= -1-x+x^5 + y(1+x+x^2+x^3)
\]
and take
\[
G(x,y)=F(x,y)+ (A+By)(1-y)\left(\frac{x^9-1}{x-1}\right)+ \ell \left(\frac{y^3-1}{y-1}\right)\left(\frac{x^9-1}{x-1}\right).
\]
This has
$G(1,1)=  3(1+9\ell)$. For the ninth roots of unity $x$ other than 1 we have $G(x,y)=F(x,y)$ 
and, writing $1$, $\om$, $\om^2$ for the cube roots of unity, we obtain
\[
\Res((1+x+x^2)(1+x^3+x^6),F(x,1)F(x,\om)F(x,\om^2))= 3^7.
\]
 This leaves the case $x=1$, $y\in\{\om,\om^2\}$, and here
\begin{equation*}
\begin{split}
&G(1,y) = -1+4y+9(1-y)(A+By)=(1-y)(-2+y+9A+9By),\\
&G(1,\om)G(1,\om^2) = 3m,
\end{split}
\end{equation*}
with
\[
m=N_1( -2+9A +(1+9B) x),\quad A,B\in \mathbb Z.
\]
Hence with $\pm G(x,y)$ we get  $\pm 3^9(1+9\ell)m$.

Note we can write
\begin{align*} 3m & =N_1( (1-x)(-2+x+9A+9Bx))\\
 & =N_1(4(1-x^2)+(-9+27A-27B)+(18B-9A)(1-x^2))\\
 & =N_1(4(1-x)+9(-1+3(A-B))+9(2B-A)(1-x)x).
\end{align*}
and $m\equiv 7\bmod 9$ is a $3$-norm of Type~1, and conversely any $3$-norm $m\equiv 7\bmod 9$ of Type~1 can be 
achieved for a suitable choice of  $A$ and $B$.

For the $m\equiv 4\bmod 9$ with a Type~1 representation, we take
\begin{equation}\label{eqnF49}
\begin{split}
F(x,y) &= (1+x+x^3+x^5+x^6) + y(1+x^3+x^5+x^8)\\
&\qquad + y^2((1+x^2)(1+x^3+x^6))
\end{split}
\end{equation}
and
\begin{equation*}\label{eqnG49}
\begin{split}
G(x,y) &= F(x,y) + x (x-1)(x^6+x^3+1)(y^2+y+1)(A + B(x^2-1))\\
&\qquad + \ell \left(\frac{y^3-1}{y-1}\right)\left(\frac{x^9-1}{x-1}\right). 
\end{split}
\end{equation*}
Here $G(1,1)=3(5+9\ell)$.
For $y$ a primitive cube root of unity $\om$ or $\om^2$, we have $G(x,y)=F(x,y)$ and 
\[
\Res(x^9-1, F(x,\om)F(x,\om^2))=3^6.
\]
For the primitive ninth roots of unity $x$, we have
\[
G(x,1)=x(1-x)(1+x^4), \quad \Res(x^6+x^3+1,G(x,1))=3.
\]
For the primitive cube roots of unity $x$, we have
\[
G(x,1)=2-5x+ x(9A(x-1)+27B) = x\left(2(1-x)+9A(x-1)+9(-1+3B)\right)
\]
and 
\[
\Res(x^2+1+1,G(x,1)) = 3m,
\]
where $m\equiv 4\bmod 9$ is a $3$-norm of Type~1.
Conversely, we can achieve any $3$-norm $m\equiv 4\bmod 9$ of Type~1 with a suitable choice of $A$ and $B$.
Hence $\pm G$ achieves all $\pm 3^9 (5+9\ell)m$.

For the $m\equiv 1\bmod 9$ with a Type~1 representation, we take
\[
F(x,y)=-1+y(1+x+x^2+x^3+x^4+x^6+x^7)
\]
and
\[
G(x,y)=F(x,y)+\left(\frac{x^9-1}{x-1}\right)y(y-1)(A+B(y^2-1))+ \ell \left(\frac{y^3-1}{y-1}\right)\left(\frac{x^9-1}{x-1}\right).
\]
For the ninth roots of unity $x\neq 1$, we have $G(x,y)=F(x,y)=-1-(x^5+x^8)y$, and for the third roots of unity $1$, $\om$, $\om^2$, we obtain
\begin{align*} &  \Res((x^6+x^3+1)(x^2+x+1), G(x,1)G(x,\om)G(x,\om^2))\\ & =\Res((x^6+x^3+1)(x^2+x+1), -1-(x^5+x^8)^3)=3^7.
\end{align*}
For $x=1$, we have $F(1,y)=-1+7y$, producing $G(1,1)=3(2+9\ell)$, and for the complex cube roots $y=\om$, $\om^2$, we have
\[
G(1,y)=-1+7y  +y\left( 9A(y-1)+27B\right)=-y\left( 1-y -9A(y-1)-9(1+3B)\right),
\]
yielding
\[
\Res(y^2+y+1,G(1,y)) = 3m,
\]
where $m\equiv 1\bmod 9$ is a $3$-norm of Type~1.
Again, by choice of $A$ and $B$ we can achieve all $3$-norms $m\equiv 1\bmod 9$ of Type~1 this way.
Hence $\pm G(x,y)$ produces $\pm 3^9(2+9\ell)m$.

It remains to deal with the $m\equiv 1\bmod 9$ with a 9-norm representation of Type~3. 
For this, we take
\[
F(x,y) = x^4+x^7+x^8+ y(1+x+x^2)
\]
and
\[
G(x,y) = F(x,y)+ (y-x^3)(x^3-1)(y-1)t(x)+  \ell \left(\frac{y^3-1}{y-1}\right)\left(\frac{x^9-1}{x-1}\right).
\]
We have $F(1,y)=3(1+y)$ giving $G(1,1)=3(2+9\ell)$ and 
\[
\Res(y^2+y+1,G(1,y)) = 3^2.
\]
For $x$ a primitive cube root of unity we have $G(x,y)=x(2+x)$, and for $y=1$, $\om$, $\om^2$ we find
\[
\Res(x^2+x+1,G(x,1)G(x,\om),G(x,\om^2)) = 3^3.
\]
If $x$ is a primitive ninth root of unity, then
\[
\Res(x^6+x^3+1, G(x,1) ) =  \Res(x^6+x^3+1, 1-x^5 )=3,
\]
and dividing the primitive cube roots of unity into those with $y=x^3$ and $y=x^6$, we have 
\[
G(x,x^3)=x(1+x)(x^2-1), \quad \Res(x^6+x^3+1, G(x,x^3)) = 3.
\]
Observing that $x^3-1=(x-1)^3u_1(x)$, $3=(x-1)^6u_2(x)$, with $u_1(x)$ and $u_2(x)$ units, we see that
\begin{align*} G(x,x^6)  & = -2x^5-x^4-x^3-2x^2-2x-1+ 3x^3(x^3-1)t(x) \\
& = (x-1)\left( 3(x^5-1)+(x^2+1)^2+ 3\left(\frac{x^3-1}{x-1}\right) x^3t(x) \right) \\
& = (1+x^2)^2(x-1) H(x),
\end{align*}
where
\[
H(x)= 1 + 3(x-1)(-x^5 - 2 x^4 +  x^2 - x - 1) + (x-1)^8 u(x)t(x),
\]
with $u(x)$ a unit.
Since $u_2(x)=-1+(x-1)t_2(x)$, we have
\[
\Res(x^6+x^3+1,G(x,x^6))=3m,
\]
with
\[
m=N_1(1 +(x-1)^7 +(x-1)^8(t_3(x)+u(x)t(x)) ).
\]
Now from \cite[Lemma 4.7]{MP23}, we have that a $9$-norm of Type~2 can be written in the form $N_2(1+(x-1)^7h(x))$,
which is a Type~3 representation when $3\nmid h(1)$.
Since for a primitive ninth root $(x^{-1}-1)^7=-(x-1)^7-(x-1)^8(x+1)$, we see that a $9$-norm of Type~3 can be written
$N_2(1+(x-1)^7+(x-1)^8h_1(x))$ for some $h_1(x)$.
Hence our $m$ has a $9$-norm representation of Type~3, and conversely any $9$-norm of Type~3 is achieved as an $m$ for a suitable choice of $t(x)$.
With $\pm G$ we achieve all $\pm 3^9(2+9\ell)m$.
By \cite[Lemmas  4.7 and 4.8]{MP23}  a $9$-norm having  a Type~3 representation must contain a Type~3 prime $p\equiv 1\bmod 9$,
so we obtain all these values by restricting to the $m=p$.

\section{All $3^{10}\parallel D$ have the required form}\label{secNoMore10}
In this and the following section, we establish that only the values described in parts~\ref{main10} and~\ref{main9} of Theorem~\ref{main} are attained as integer group determinants.
We write
\[
F(x,y) = \sum_{0\leq j\leq 8} a_j(x-1)^j + (y-1) \sum_{0\leq j \leq 8} b_j(x-1)^j + (y-1)^2\sum_{0\leq j\leq 8}c_j (x-1)^j,
\]
and split the determinant into the product  of $F(1,1)$, the four $3$-norms
\[
A_1=N_1(F(x,1)),\; A_2=N_1(F(x,x)),\; A_3=N_1(F(x,x^2)), \; A_4=N_1(F(1,y)),
\]
and the three $9$-norms
\[
B_1=N_2(F(x,1)),\; B_2=N_2(F(x,x^3)),\; B_3=N_2(F(x,x^6)).
\]
Suppose that $3^{9}\parallel D$ or $3^{10}\parallel D$.
Since $F(x,y)\equiv F(1,1) \bmod {\langle x-1,y-1 \rangle }$, we have
\[
3\mid D \implies 3\mid F(1,1),A_1,A_2,A_3,A_4,B_1,B_2,B_3.
\]
Writing 
\[
F(x,y)=a_0+a_1(x-1)+b_0(y-1)  \bmod  {\langle (x-1)^2, (y-1)^2, (x-1)(y-1)  \rangle},
\]
observe that if $3\mid a_1$, then for the primitive third and ninth roots of unity $x$ we have
\[
(x-1)^2\mid F(x,1) \implies 3^2\mid A_1,B_1,
\]
while for the primitive ninth roots, $(x-1)^3\mid (x^3-1)$ and 
\[
(x-1)^2 \mid F(x,x^3),F(x,x^6) \implies 3^2 \mid B_2,B_3,
\]
so $3^{12}\mid D$.
We can therefore assume that $3\nmid a_0$.
This ensures that
$(x-1)\parallel F(x,1)$ for the primitive third and ninth roots, $(x-1)\parallel F(x,x^3), F(x,x^6)$ for the primitive ninth roots, and
\[
3 \parallel A_1, B_1, B_2, B_3.
\]
Notice that if $3\mid b_0$,  then for $x$ or $y$ a primitive cube root of unity, we have
\[
3\mid b_0 \implies (x-1)\parallel F(x,x),F(x,x^2), \; (y-1)^2\mid F(1,y) \implies 3\parallel A_2,A_3, \; 3^2\mid A_4,
\]
and, since $b_0 (x^2-1) = 2b_0(x-1)+b_0(x-1)^2$, we find
\[
b_0\equiv -a_0\bmod 3 \implies (x-1)^2 | F(x,x), \; (x-1)\parallel F(x,x^2), \; (y-1)\parallel F(1,y),
\]
giving $3^2\mid A_2,\; 3\parallel A_3,A_4$, while
\[
b_0\equiv a_0\bmod 3 \implies (x-1)^2 \mid F(x,x^2), \; (x-1)\parallel F(x,x), \; (y-1)\parallel F(1,y),
\]
yielding $3^2\mid A_3,\; 3\parallel A_2,A_4$.
By replacing $F(x,y)$ by $F(xy,y)$ or $F(xy^2,y) \bmod {y^3-1}$, we run through the same values, but now we may assume that 
\[
3^2\mid A_4, \quad 3\parallel A_1,A_2,A_3,B_1,B_2,B_3,
\]
and by replacing $F$ by $-F$ if necessary, we may assume that
\[
a_1 \equiv 1\bmod 3,\;\; 3\mid b_0 \implies a_1 = 1+3\alpha_1, \: b_0 = 3\beta_0.
\]
Thus, for $3^{10}\parallel D$, we have one of the following two cases.
\begin{itemize}
\item Case~1: $3^2\parallel F(1,1)=a_0$, $3^2\parallel A_4$, and
$(y-1)^2  \parallel F(1,y)=a_0+3\beta_0(y-1)+c_0(y-1)^2$, which implies $3\nmid c_0$.
\item Case~2: $3\parallel F(1,1)=a_0=3\alpha_0$, $3\nmid \alpha_0$, $3^3\parallel A_4$, and
$(y-1)^3\parallel  F(1,y)=3\alpha_0+3\beta_0(y-1)-3c_0y$, which implies $\alpha_0\equiv c_0 \bmod 3$, $\beta_0\not\equiv c_0 \bmod 3$.
\end{itemize}
In both cases, we have that $c_0\equiv \pm1 \bmod 3$.

Suppose that we had a determinant $D=3^{10}m$, $3\nmid m$, with $m$ containing no Type~1 prime $p\equiv 1\bmod 9$.
Then the $9$-norms $B_1$, $B_2$, $B_3$ cannot be written $B_i=3N_2(\alpha_i)$ with $\alpha_i$ of Type~1.
That is, for a primitive ninth root of unity $x$, if we write
\begin{equation*}
\begin{split}
F(x,1) &= (x-1)\left(F_1(x) + (x-1)^6 g_1(x)\right),\\
F(x,x^3) &= (x-1)\left(F_2(x) + (x-1)^6 g_2(x)\right),\\
F(x,x^6) &= (x-1)\left(F_3(x) + (x-1)^6 g_3(x)\right),
\end{split}
\end{equation*}
then the $F_i(x)$ have only Type~2 representations (likewise their products).

We reduce $F(x,1)$, $F(x,x^3)$ and $F(x,x^6) \bmod (x-1)^7$.
Since 
\begin{equation*}
\begin{split}
3 &= (x-1)(x^2-1)(x^4-1)(x^5-1)(x^7-1)(x^8-1)\\
&= (x-1)^6 g_1(x) = -(x-1)^6+(x-1)^7g_2(x),
\end{split}
\end{equation*}
we can reduce the coefficients mod $3$ for any terms with at least one $(x-1)$.
We can write $x^3-1 = (x-1)^3+3x(x-1)$, and $x^6-1 = -(x^3-1)+(x^3-1)^2+3(x^3-1)$.
Hence with
\be \label{def0} F_0 = 1+a_2(x-1) +a_3(x-1)^2+a_4(x-1)^3+a_5(x-1)^4+a_6(x-1)^5, \ee
in Case~2 we may write $a_0=3\alpha_0= -\alpha_0(x-1)^6 \bmod{ (x-1)^7}$, and we have
\begin{equation}\label{only3form} 
\begin{split}
F_1 &= F_0-\alpha_0(x-1)^5,\\
F_2 &= F_1+ b_1(x-1)^3+b_2(x-1)^4+b_3(x-1)^5 + c_0(x-1)^5,\\
F_3 &= F_1- b_1(x-1)^3- b_2(x-1)^4- b_3(x-1)^5 + c_0(x-1)^5,
\end{split}
\end{equation}
with 
\begin{equation*}\label{range} 
a_2,a_3,a_4,a_5,a_6 ,b_1,b_2,b_3\in \{0,-1,1\}, \quad c_0 = \pm 1.
\end{equation*}
Case~1 is similar, except that $a_0\equiv 0 \bmod{ (x-1)^{12}}$, and we take $F_1=F_0$.

Observe that
\[
F_2F_3  \equiv \left(F_1+c_0(x-1)^5\right)^2 \equiv F_1^2+2c_0(x-1)^5 \bmod (x-1)^6.
\]
As $F_1$  has no Type~1 representation, we can write $u(x)F_1=1+(x-1)^6g_1(x)$ for some unit $u(x)$, but then 
\[
u(x)^2F_2F_3=1+2c_0u(1)^2(x-1)^5+(x-1)^6g_3(x), \quad 3\nmid 2c_0u(1)^2,
\]
gives a Type~1 representation.
Hence at least one of $F_2$, $F_3$ must have a Type~1 representation, producing a contradiction.

\section{All  $3^9\parallel D$ have the required form}\label{secNoMore9}

Suppose $3^9\parallel D$. We proceed as before, assuming $3^2\mid A_4$.
This time we deduce that $3\parallel F(1,1)$ and $3^2\parallel A_4$.
That is, $a_0=3\alpha_0$, $3 \nmid \alpha_0$, $b_0=3\beta_0$, and
\begin{gather*}
B_1 = 3m_1,\;
B_2 = 3m_2,\;
B_3 = 3m_3,\;
3\nmid m_1m_2m_3,\\
A_1 = 3n_1,\;
A_2 = 3n_2,\;
A_3 = 3n_3,\;
3\nmid n_1n_2n_3,\\
A_4 = 9n_0,\;
n_0 = N_1(\alpha_0 + \beta_0(y-1) - c_0y),\;
\alpha_0\not\equiv c_0 \bmod 3.
\end{gather*}
Suppose that $D=3^9m$ is not of  the form \ref{main9i}, \ref{main9ii} or \ref{main9iii} in Theorem~\ref{main}, so $m\equiv \pm 7\bmod 9$.
Since the $9$-norms $m_i$ consist of $q^f\equiv 1\bmod 9$, we have that the $m_i\equiv 1\bmod 9$, and the $m\equiv \pm 7\bmod 9$ reduces to
\be \label{7cong} n_0n_1n_2n_3 \equiv \pm 2 \bmod 9. \ee
Moreover, $m$ does not contain a Type~1 prime $p\equiv1$, $4$ or $7\bmod 9$, or a Type~3 prime $p\equiv 1\bmod 9$.
That tells us that the three $9$-norms have only Type~2 representations, 
\be \label{step1} m_i=N_2(1+(x-1)^5h_i(x)) \implies 3\mid h_i(1), \ee
and moreover the Type~2 representations must be Type~4,
\be \label{step3} m_i= N_2( 1+(x-1)^7h_i(x)) \implies 3\mid h_i(1).\ee
Similarly the $3$-norms have only Type~2 representations
\begin{equation*}\label{step2} 
n_i = N_1( A+3B(x-1)) \implies 3\mid B.
\end{equation*}

\subsection{When all the 9-norms are Type~4}

As in the previous section, to apply \eqref{step1} we write
\begin{align*}
F(x,1) &= (x-1)F_1(x) \bmod (x-1)^7,\\
F(x,x^3) &= (x-1)F_2(x) \bmod (x-1)^7,\\
F(x,x^6) &= (x-1)F_3(x) \bmod (x-1)^7,
\end{align*}
and we obtain a contradiction if we can show that one of $F_1(x)$, $F_2(x)$, $F_3(x)$ has a Type~1 representation.  
We again have \eqref{only3form}, and the same argument produces a contradiction if $3\nmid c_0$.
Hence we can assume that $3\mid c_0$ and write $c_0=3\gamma_0$, and take
\begin{align*}
F_1 &= F_0 - \alpha_0(x-1)^5,\\
F_2 &= F_1 + b_1(x-1)^3 + b_2(x-1)^4 + b_3(x-1)^5,\\
F_3 &= F_1 - b_1(x-1)^3 - b_2(x-1)^4 - b_3(x-1)^5,
\end{align*}
with $F_0$ defined as in \eqref{def0}.
Note that
\[
F_1 F_2F_3 = F_1^3 = 1+a_1^3(x-1)^3 \bmod (x-1)^6,
\]
and we find $(\ve_1,\ldots ,\ve_5)=(0,0,a_1,0,0)$ in the database described in Section~\ref{secAlg} is Type~2 for $a_1=0$, $-1$ or $1$.
In particular, if $F_1$  and $F_2$ are both Type~2, then $F_3$ is automatically also Type~2.

Writing in the form \eqref{standardform}, we get
\begin{equation} \label{Type1Test}
\begin{split}
F_1 &=1+a_2\pi +a_3\pi^2+a_4\pi^3+a_5\pi^4+(a_6-\alpha_0 \text{ mod } 3)\pi^5 \bmod \pi^6,\\ 
F_2 &= 1+a_2 \pi + a_3\pi^2 + (a_4+b_1 \text{ mod }3)\pi^3+(a_5+b_2 \text{ mod }3)\pi^4\\
&\qquad + (a_6+b_3-\alpha_0\text{ mod } 3)\pi^5 \bmod \pi^6, 
\end{split}
\end{equation}
where $\alpha_0$, $a_2$, $a_3$, $a_4$, $a_5$, $a_6$, $b_1$, $b_2$ are defined mod $3$ so can be assumed to be in $\{0,-1,1\}$, with $a_0=\pm 1$.

To apply \eqref{step3} for Type~3 primes, we write
\begin{align*}
F(x,1) &= (x-1)H_1(x) \bmod (x-1)^9,\\
F(x,x^3) &= (x-1)H_2(x) \bmod (x-1)^9,\\
F(x,x^6) &= (x-1)H_3(x) \bmod (x-1)^9.
\end{align*}
We can take
\begin{align*}
H_1 &= -\alpha_0 (x-1)^5+ (1+3\alpha_1)+a_2(x-1)+a_3(x-1)^2+\cdots +a_8(x-1)^7, \\
H_2 &= H_1 + (x^3-1)(b_1+b_2(x-1)+b_3(x-1)^2+b_4(x-1)^3 + b_5(x-1)^4)\\
&\qquad + (x^3-1)^2 ( c_1  + c_2(x-1)),\\
 H_3 &=  H_1 + (x^6-1)(b_1+b_2(x-1)+b_3(x-1)^2+b_4(x-1)^3 + b_5(x-1)^4)\\
 &\qquad + (x^3-1)^2 (c_1+c_2(x-1)).
\end{align*}
Here $\alpha_0,\alpha_1,a_3,\ldots,a_8,b_1\ldots,b_5,c_1,c_2$ are defined mod $3$, and $a_2=e_2+3e_3$ are defined mod $9$.
With 
\begin{align*}
x^3-1 &= \pi^3-\pi^7 \bmod \pi^8,\\
x^6-1 &= -\pi^3+\pi^6+\pi^7 \bmod \pi^8,
\end{align*}
writing in powers of $\pi$ we get
\begin{equation}\label{Step3Test}
\begin{split}
H_1 &= 1+e_2\pi +a_3\pi^2+a_4\pi^3+a_5\pi^4+(a_6-\alpha_0 \bmod 3)\pi^5\\
&\qquad +  (a_7-\alpha_1 \bmod 3)\pi^6 +(a_8-e_3 \bmod 3)\pi^7 \bmod \pi^8,\\
H_2 &= 1+e_2 \pi + a_3\pi^2 + (a_4+b_1 \bmod 3)\pi^3+(a_5+b_2 \bmod 3)\pi^4\\
&\qquad + (a_6+b_3-\alpha_0 \bmod 3)\pi^5 +  (a_7-\alpha_1+b_4+c_1 \bmod 3)\pi^6\\
&\qquad +(a_8-e_3+b_5-b_1+c_2 \bmod 3)\pi^7 \bmod \pi^8,\\
H_3 &= 1+e_2 \pi + a_3\pi^2 + (a_4-b_1 \bmod 3)\pi^3+(a_5-b_2 \bmod 3)\pi^4\\
&\qquad + (a_6-b_3-\alpha_0 \bmod 3)\pi^5 + (a_7-\alpha_1-b_4+c_1+b_1 \bmod 3)\pi^6\\
&\qquad +(a_8-e_3-b_5+b_2+b_1+c_2 \bmod 3)\pi^7 \bmod \pi^8.
\end{split}
\end{equation}

\subsection{When all the 3-norms are Type~2}

This time, with $x$ a primitive cube root of unity, we write 
\begin{equation*}
\begin{split}
F(x,1) &\equiv (x-1)G_1  \bmod (x-1)^5,\\
F(x,x) &\equiv (x-1)G_2 \bmod (x-1)^5,\\
F(x,x^2) &\equiv (x-1)G_3 \bmod (x-1)^5,
\end{split}
\end{equation*}
we write the $G_i$ in the form \eqref{Type1form} mod $(x-1)^4$, and we apply \eqref{Step2} to each of them.

With $\alpha_0$, $a_2$, $b_1$ determined mod $9$, and $a_3$, $a_4$, $b_2$, $b_3$, $c_1$, $c_2$, $\beta_0$, $\gamma_0$ determined mod $3$, we have
\begin{align*}
G_1 &= \alpha_0 (x^2-1)+(1+3\alpha_1)+a_2(x-1)+a_3(x-1)^2+a_4(x-1)^3,\\
G_2 &= G_1+(3\beta_0+b_1(x-1)+b_2(x-1)^2+b_3(x-1)^3)\\
&\qquad + (x-1)(3\gamma_0+c_1(x-1)+c_2(x-1)^2),\\
G_2 &= G_1+(x+1)(3\beta_0+b_1(x-1)+b_2(x-1)^2+b_3(x-1)^3) \\
&\qquad +(x-1)(x+1)^2(3\gamma_0+c_1(x-1)+c_2(x-1)^2).
\end{align*}
With $3=O((x-1)^2)$ and writing 
\begin{align*}
(x^2-1) &= -(x-1)-3,\\
(x-1)^2 &= -3-3(x-1),\\
(x-1)^3 &= -3(x-1) -3(x-1)^2,
\end{align*}
we have
\[
G_1 = (1-3\alpha_0+3\alpha_1-3a_3)+(x-1)(a_2 -\alpha_0-3a_3-3a_4) \bmod (x-1)^4,
\]
and \eqref{Step2} becomes
\begin{equation}\label{defininga2}
\begin{split}
a_2 \equiv \alpha_0 \bmod 3 & \implies a_2 \equiv \alpha_0 +3a_3+3a_4 \bmod 9,\\
a_2 \equiv \alpha_0 +1 \bmod 3 & \implies a_2 \equiv 1-2\alpha_0 +3\alpha_1 +3a_4 \bmod 9,\\
a_2 \equiv \alpha_0 -1 \bmod 3 & \implies a_2 \equiv 4\alpha_0 -4 -3\alpha_1-3a_3 +3a_4 \bmod 9,
\end{split}
\end{equation}
while
\begin{align*}
G_2 &= (1-3\alpha_0+3\alpha_1-3a_3+3\beta_0-3b_2-3c_1)\\
&\quad + (x-1)(a_2-\alpha_0 +b_1-3a_3-3a_4-3b_2-3b_3+3\gamma_0-3c_1-3c_2) \bmod {(x-1)^4}, 
\end{align*}
and \eqref{Step2} becomes
\begin{equation}\label{definingb1}
\begin{split}
b_1 \equiv \alpha_0 -a_2 \bmod 3  \implies b_1 &\equiv \alpha_0-a_2 +3a_3+3a_4 +3b_2+3b_3\\  & \quad  -3\gamma_0+3c_1 +3c_2\bmod 9,\\
b_1 \equiv \alpha_0 -a_2 +1 \bmod 3  \implies b_1 &\equiv1-2\alpha_0 -a_2 +3a_4 + 3\alpha_1 +3\beta_0+3b_3\\   & \quad -3\gamma_0+3c_2 \bmod 9,\\
b_1 \equiv \alpha_0 -a_2 -1 \bmod 3  \implies b_1 &\equiv 4\alpha_0 -4-a_2  -3\alpha_1  -3a_3 +3a_4 -3\beta_0\\ & \quad -3b_2 +3b_3-3\gamma_0-3c_1+3c_2\bmod 9.
\end{split}
\end{equation}
Using
\begin{gather*}
x+1=2+(x-1), \quad (x+1)(x-1)^2=3, \quad (x+1)(x-1)^3=3(x-1),\\
(x+1)^2(x-1)=-2(x-1)-3,  \quad (x+1)^2(x-1)^2=6+3(x-1),\\
(x+1)^2(x-1)^3=-3(x-1)-9,
\end{gather*}
we obtain
\begin{align*}
G_3 &= (1-3\alpha_0+3\alpha_1-3a_3-3\beta_0-3b_1+3b_2-3c_1)\\
&\quad +(x-1)(a_2-\alpha_0-b_1-3a_3-3a_4+3\beta_0+3b_3+3\gamma_0+3c_1-3c_2),
\end{align*}
and \eqref{Step2} becomes
\begin{equation}\label{checkcong}
\begin{split}
b_1 \equiv a_2-\alpha_0 \bmod 3  \implies b_1 &\equiv a_2-\alpha_0 -3a_3-3a_4 +3\beta_0+3b_3 \\ & \quad +3\gamma_0+3c_1 -3c_2\bmod 9,\\
b_1 \equiv a_2-\alpha_0 -1 \bmod 3  \implies b_1 &\equiv 4a_2-4-\alpha_0-3a_4 - 3\alpha_1 -3\beta_0-3b_2\\  & \quad +3b_3+3\gamma_0-3c_1-3c_2 \bmod 9,\\
b_1 \equiv a_2-\alpha_0 +1 \bmod 3  \implies b_1 &\equiv 1-\alpha_0 -2a_2  +3\alpha_1  +3a_3 -3a_4 \\ &\quad +3b_2+3b_3+3\gamma_0-3c_2\bmod 9.
\end{split}
\end{equation}
Notice $\alpha$ mod $ (x-1)^3$ determines $N_1(\alpha) \bmod 9$, so
$n_i\equiv N_1(g_i)$ mod 9 with
\begin{align*}
g_1 &= 1-3\alpha_0+3\alpha_1-3a_3 +(x-1)(a_2-\alpha_0),\\
g_2 &= g_1-3c_1 +3(\beta_0-b_2)+b_1(x-1),\\
g_2 &= g_1-3c_1  -3(\beta_0-b_2) -b_1(x-1)-3b_1,
\end{align*}
yielding
\begin{align*}
g_1g_2g_3 &\equiv g_1 \left((g_1-3c_1)^2-(3(\beta_0-b_2)+b_1(x-1))^2 -3b_1g_1\right)\\ 
&= g_1\left( g_1^2+3c_1+3b_1^2-3b_1 \right) \bmod (x-1)^3 \\
&= 1+3c_1+3b_1^2 -3b_1 \bmod (x-1)^3,
\end{align*}
and 
\[
n_1n_2n_3 \equiv 1  -3c_1-3b_1^2+3b_1 \bmod 9.
\]
Writing $\alpha_0=e_0+3e_1$ with $e_0=\pm 1$, $e_1=0,\pm 1$, we have 
\begin{align*}
\alpha_0n_0 &= \alpha_0 N_1( \alpha_0-3\gamma_0+(\beta_0-3\gamma_0)(y-1))\\
&= \alpha_0 \left( (\alpha_0-3\gamma_0)^2+3(\beta_0-3\gamma_0)^2 -3(\alpha_0-3\gamma_0)(\beta_0-3\gamma_0)\right)\\
&= \alpha_0\left( \alpha_0^2+3\alpha_0\gamma_0+3\beta_0^2-3\alpha_0\beta_0\right)\\
&\equiv e_0+3\gamma_0+3\beta_0^2e_0-3\beta_0 \bmod 9.
\end{align*}
Thus \eqref{7cong} becomes
\[
e_0 -e_0(3c_1+3b_1^2-3b_1) + 3\gamma_0 +3\beta_0^2e_0-3\beta _0 \equiv \pm 2 =-2e_0 \bmod 9,
\]
giving
\be \label{definingc1} c_1\equiv 1-b_1^2+b_1 +\beta_0^2 +(\gamma_0-\beta_0)e_0 \bmod 3. \ee

\subsection{Completing the proof of Theorem~\ref{main}}\label{subsecCompletingProof}

The proof is computational.
We break the process into three steps.

\subsubsection*{Step 1}
We select $\alpha_0=\pm 1$ and $a_2$, $a_3$, $a_4$, $a_5$, $a_6$, $b_1$, $b_2$, $b_3$ from $\{0,-1,1\}$
and check the sequences $(\ve_1,\ldots,\ve_5)$ for $F_1$ and $F_2$ in \eqref{Type1Test} against the Type~1 database described in Section~\ref{secAlg}.
If either are Type~1, we discard that selection.
If $F_1$ and $F_2$ are both Type~2 then this assignment of values proceeds to Step~2.
Of the $2\cdot 3^8=13122$ choices, all but $1458$ are eliminated.

\subsubsection*{Step 2}
We take the mod $3$ values $\alpha_0$, $a_2$, $a_3$, $a_4$, $b_1$, $b_2$, $b_3$ in $\{0,-1,1\}$ that passed Step~1 as input.
We take $\alpha_0=e_0+3e_1 \bmod 9$ with $e_0$ the $\alpha_0 \bmod 3$ value from Step~1.
For the $3^5$ choices of $e_1$, $\alpha_1$, $\beta_0$, $\gamma_0$, $c_2$ from $\{0,-1,1\}$, we compute $c_1$ in $\{0,-1,1\}$ using \eqref{definingc1}, we compute $a_2\bmod 9$ using \eqref{defininga2}, and $b_1\bmod 9$ using \eqref{definingb1}.
We then check whether \eqref{checkcong} holds; if not we discard that assignment of values, if it holds we go to Step~3.
Of the $3^5\cdot1458 = 354294$ possibilities, exactly $2/3$ were eliminated.

\subsubsection*{Step 3}
We take the mod $3$ values $a_3$, \ldots, $a_6$, $b_1$, $b_2$, $b_3$, $\alpha_1$, $c_1$, $c_2$ from $\{0,-1,1\}$ and $\alpha_0=\pm 1$ from Steps~1 and~2, and the $a_2=e_2+3e_3 \bmod 9$ from Step~2.
For each of the $3^4$ values of $a_7$, $a_8$, $b_4$, $b_5$ from $\{0,1,-1\}$, we check the sequences $(\ve_1,\ldots ,\ve_7)$ for $G_1$, $G_2$, $G_3$ from \eqref{Step3Test} against the Type~3 database.
All $118098$ values that survived Step~2 were eliminated in Step~3.
Of these, $78732$ were eliminated by $H_1$, which required testing $9\cdot78732=708588$ sequences, since $H_1$ does not depend on $b_4$ or $b_5$.
The other 39366 were eliminated by $81\cdot39366=3188646$ tests using $H_2$ and $H_3$.
The appeal to $H_2$ eliminated $2/3$ of these, and the check using $H_3$ (which was performed only if the $H_2$ test failed) eliminated the other $1/3$.

\subsection*{Data Availability}
A record of the computation described in Section~\ref{subsecCompletingProof} is available at \cite{GithubData}, with the outcomes of all tests performed in Steps~1, 2 and~3.
Lists of primes of Types~1, 2, 3 and~4 up to $50000$ are available there too, as well as the lists over $\{-1,0,1\}$ that are used in Section~\ref{subsec9norms} for classifying certain algebraic integers as one of these types.

\subsection*{AI Statement}
GPT-5.6 Sol was used to determine each of the polynomials $f_1(x)$, $f_2(x)$ and $f_3(x)$ in Section~\ref{subsecAchieve10}, and GPT-6 Astra for the polynomial \eqref{eqnF49} in Section~\ref{subsecAchieve9}.

\bibliographystyle{amsplain}

\end{document}